\documentclass[11pt]{article}
\usepackage{amssymb}
\parskip \baselineskip
\newcommand{\ds}{\displaystyle}

\newcommand{\Z}{\mathbb{Z}}

\newcommand{\ol}{\overline}

\newcommand{\ra}{\rightarrow}
\newcommand{\Ra}{\Rightarrow}

\begin{document}

\title{How to partition or count an abstract simplicial complex, given its facets}

\author{Marcel Wild}

\maketitle

\begin{quote}
{\bf Abstract}. Given are the facets of an abstract (finite) simplicial complex ${\cal S}{\cal C}$. We show how to partition ${\cal S}{\cal C}$ into few pieces, each one compactly encoded by the use of wildcards. Such a representation is useful for the optimization of target functions $\ {\cal S}{\cal C} \ra \Z$. Other applications concern combinatorial commutative algebra, the speed up of inclusion-exclusion, and Frequent Set Mining. Merely calculating the face-numbers of ${\cal S}{\cal C}$ can be done faster than partitioning ${\cal S}{\cal C}$. Our method compares favorably to inclusion-exclusion and binary decision diagrams.
 \end{quote}

\section{Introduction}

A {\it simplicial complex} (also called {\it set ideal}) based on a set $W$ of cardinality $w$ is a family ${\cal S}{\cal C}$ of subsets $X \subseteq W$ (called {\it faces}) such that from $X \in {\cal S}{\cal C}$, $Y \subseteq X$, follows $Y \in {\cal S}{\cal C}$. In this article all structures are assumed to be finite. In particular, all simplicial complexes ${\cal S}{\cal C}$ contain maximal faces, called the {\it facets} $F_i$ of ${\cal S}{\cal C}\, (1\leq i \leq h)$. The facets uniquely determine ${\cal S}{\cal C}$.  

This article is about partitioning ${\cal S}{\cal C}$ into $R$ pieces $r_i \subseteq {\cal S}{\cal C}$. Usually $R$ is small compared to $|{\cal S}{\cal C}|$ and by the use of wildcards each $r_i$ packs its many faces in a compact way. Roughly speaking $r_i$ is a {\it multivalued} {\it row} of length $w$ with components $0,1,2,e$. Sections 3 to 5 present methods geared to arbitrary simplicial complexes given by their facets. Sections 6 to 8 give applications of it all to inclusion-exclusion, to combinatorial commutative algebra and to Frequent Set Mining.
More specifically, the section break up is as follows. 

In Section 2 the facets are used, in junction with binary decision diagrams, to calculate $N=|{\cal S}{\cal C}|$.

In Section 3 multivalued rows enter the picture in order to refine the calculation of $N$ to the calculation of {\it all} numbers $f_k$ of $k$-element faces $(1\leq k \leq w)$. This entails applying a certain $e$-algorithm to $h$ constraints coupled to the $h$ facets. In 3.1 we show that this approach 
outperforms inclusion-exclusion and binary decision diagrams. In 3.2 we briefly discuss the complexity of calculating $f_k$ for a {\it specific} value $k \in [w]$.

In Section 4 we invest $\frac{1}{2}(h^2-h)$ constraints in order to represent ${\cal S}{\cal C}$ {\it itself} as a disjoint union of few multivalued rows (which yields all $f_k$'s as a side product). 

Section 5 shows how to get ${\cal S}{\cal C}$ partitioned in case not the facets but the {\it minimal non-faces} of ${\cal S}{\cal C}$ are provided.  In this context a dual version of the $e$-algorithm, called $n$-algorithm applies.

Section 6 indicates how the $n$-algorithm can be used to skip a potentially large number of zero terms in the familiar length $2^n$ inclusion-exclusion expansion.

In Section 7 the $e$-algorithm is used to calculate the following objects that come along with the kind of simplicial complexes appearing in combinatorial commutative algebra: $h$-polynomials, reduced homology groups, and links of faces.

In Section 8 simplicial complexes coming from Frequent Set Mining are sliced in two ways: cardinality-wise and frequency-wise. This allows to calculate certain relevant probabilities. Furthermore we touch upon {\it closed} frequent sets and give a quick proof of the main result in [UAUA].

\section{Calculating the cardinality of a simplicial complex from its facets}

Given the facets of an arbitrary simplicial complex ${\cal S}{\cal C}$, how can we compute its cardinality? That is of considerable interest e.g. for Frequent Set Mining (Sec.8). 
Say $W = [14]:= \{1,2, \cdots, 14\}$, and ${\cal S}{\cal C}'$ is defined by its facets

$F_1 = \{ 1, 2, 5, 6, 7, 8, 10, 11, 12, 13, 14 \}, \quad F_2 = \{1, 2, 3, 4, 6, 7, 8, 9, 11, 12, 13, 14 \}$

$F_3 = \{1, 2, 3, 4, 5, 8,  9, 10, 13, 14 \}, \quad F_4 = \{1,2, 3, 4, 5, 6, 7, 9, 10, 11, 12 \}$

$F_5 = \{9, 10, 11, 12, 13, 14 \}, \quad F_6 = \{1, 2, 6, 7, 9, 10, 11, 14 \}$.

Writing ${\cal P}(X)$ for the powerset of $X$ inclusion-exclusion yields

(1) \quad $|{\cal S}{\cal C}'| = |{\cal P}(F_1) \cup \cdots \cup {\cal P}(F_6)| = \ds\sum_{i=1}^6 |{\cal P}(F_i)| - \ds\sum \{|{\cal P}(F_i) \cap {\cal P}(F_j)|: 1 \leq i < j \leq 6\}$

\hspace*{1cm} $+ \cdots - \cdots + |{\cal P}(F_1) \cap \cdots \cap {\cal P}(F_6)|$.

Since say $|{\cal P}(F_1) \cap {\cal P}(F_2)| = |{\cal P}(F_1 \cap F_2)|=2^9$, adding and substracting $64$ powers of 2 yields $|{\cal S}{\cal C}'| = 7600$. Unfortunately the inclusion-exclusion approach costs $O(2^hw)$ time units (in the usual RAM model). This is feasible only for small numbers $h$ of facets.

Consider instead the Boolean function $b: \ \{0,1\}^W \ra \{0,1\}$ whose $i$th conjunction consists of the negated literals with indices {\it not} in the $i$th facet:

$$b(x) \quad =\quad (\ol{x}_3 \wedge \ol{x}_4 \wedge \ol{x}_9) \vee (\ol{x}_5 \wedge \ol{x}_{10}) \vee \cdots \vee (\ol{x}_3 \wedge \ol{x}_4 \wedge \ol{x}_5 \wedge \ol{x}_8 \wedge \ol{x}_{12} \wedge \ol{x}_{13}).$$
It follows that $b(x) =1$ if and only if the {\it support} of $x$, i.e. $X := \{i|x_i = 1\}$, belongs to ${\cal S}{\cal C}'$. Therefore 
$$|{\cal S}{\cal C}'| = {\tt SatisfiabilityCount}[b] = 7600,$$

where {\tt SatisfiabilityCount} is a Mathematica-command that is based on the technique of binary decision diagrams (BDD), see [K]. Generally, if a counting problem that doesn't immediately yield to combinatorial evaluation can be recast as counting the number of models of a Boolean function, then BDD's are likely the fastest option. However, it is hard to assess the performance of BDD's on a theoretic level.

\section{Getting the face numbers of a simplicial complex from its facets}

How to compute the face numbers $f_k$ of a simplicial complex that is given by its facets? Inclusion-exclusion still works, provided that in (1) instead of ${\cal P}(F_i)$ we use the $k$-slice ${\cal P}(F_i,k) : = \{X \subseteq F_i: \ |X| = k\}$ throughout. Then the powers of 2 in (1) become some obvious binomial coefficients. But inclusion-exclusion is as inefficient as before. Unfortunately also BDD's are ruled out because incorporating a cardinality constraint into a Boolean function blows it out of proportion (see [W1]).
 
However, one can proceed as follows. Returning to our example, consider the complements $H_i: = W \setminus F_i$ of the facets $(1\leq i \leq 6)$ and observe that for all $X \in {\cal P}(W)$ one has:
$$X \not\in {\cal S}{\cal C}' \quad \Leftrightarrow \quad (\forall i) (X \not\subseteq F_i) \quad \Leftrightarrow  \quad (\forall i)( X \cap H_i \neq \emptyset).$$
Thus the complementary set filter ${\cal S} {\cal F}' := {\cal P}(W)\setminus {\cal S}{\cal C}'$ can be viewed as  the {\it transversal hypergraph} ${\cal T}r({\cal H}')$, i.e. as the family of transversals of the hypergraph ${\cal H}' = \{H_1, \cdots, H_6\}$. 
Because each $k$-element subset of $W$ is in exactly one of ${\cal S}{\cal C}'$ and ${\cal S}{\cal F}'$ we deduce that

(2) \quad $f_k = {w\choose k} - \tau_k \quad (1 \leq k \leq w),$

where 
\begin{center}
$\tau_k : =$  number of $k$-element transversals of the hypergraph at hand (here ${\cal H}')$.
\end{center}
The point is that ${\cal S}{\cal F}' = {\cal T}r({\cal H}')$ is representable as a disjoint union of so called $\{0,1,2,e\}$-{\it valued}\footnote{Interchangeably we speak of {\it multivalued rows} or just {\it rows}.} rows $r_1$ to $r_7$ (Table 1) from which the numbers $\tau_k$ are easily computable. We only discuss {\it what} this representation achieves, and refer to [W1] for {\it how} it is obtained.

\begin{tabular}{c|c|c|c|c|c|c|c|c|c|c|c|c|c|c|} 
& 1 & 2 & 3 & 4& 5 &6 &7 & 8 & 9 & 10 & 11 & 12 & 13 & 14 \\ \hline
$r_1=$ & 2 & 2 & $e_1$ & $e_1$ & $e_2$ & $e_3$ & $e_3$ & $e_4$ & 2 & $e_2$ & $e_3$ & $e_3$ & $e_4$ & $e_4$ \\ \hline
$r_2=$ & 2 & 2 & 0 & 0 & 1 & $e_3$ & $e_3$ & $e_4$ & 1 & 2& $e_3$ & $e_3$ & $e_4$ & $e_4$\\ \hline
$r_3=$ & 2 & 2 & 0 & 0 & 0 & $e_1$ & $e_1$ & $e_2$ & 1 & 1 & 2 & 2 & $e_2$ & 2 \\ \hline
$r_4=$ & 2 & 2 &0 & 0 & 0 & $e$ & $e$ & 0 & 1 & 1 & 2 & 1 & 0 & 1 \\ \hline
$r_5=$ & 2 & 2 & 0 & 0 & 0& 0 & 0 & 1 & 1 & 1 & $e$ & $e$ & 2 & 2\\ \hline
$r_6=$ & $e_1$ & $e_1$ & 0 & 0 & 0 & 0 & 0 & 0 & 1 & 1& 2 & 1 & $e_2$ & $e_2$ \\ \hline
$r_7=$ & $e$ & $e$ & 0 & 0 & 0 & 0 & 0 & 0 & 1 & 1 & 1 & 0 & 1 & 2 \\ \hline
\end{tabular}

Table 1

Namely, each $r_i$ comprises a bunch of $0,1$-strings $x$ whose supports $X$ are transversals of ${\cal H}'$. Besides the {\it don't care} symbol 2, which can be freely chosen to be either $0$ or $1$, we use the wildcard $e e \cdots e$ which means ``at least one 1''.  In other words, only $00 \cdots 0$ is forbidden.  If several such $e${\it -bubbles} appear within a row, they are distinguished by indices. Thus say $x=(1,0,1,1,1,1,0,0,0,0,0,0,1,1)$ is a member of $r_1$, generated by $e_1 e_1 =11, e_2e_2=10$, $e_3e_3e_3e_3=1000$, and $e_4e_4e_4=011$. Correspondingly  $X = \{1,3,4,5,6,13, 14\} \in {\cal T}r({\cal H}')$.

Here is the general (up to permutation of the entries) definition of a $\{0,1,2,e\}${\it -valued row}:

(3) \quad $r\quad =\quad (\underbrace{0, \cdots, 0}_\alpha , \underbrace{1, \cdots, 1}_\beta , \underbrace{2, \cdots, 2}_\gamma , \underbrace{e_1, \cdots, e_1}_{\varepsilon_1} , \cdots , \underbrace{e_t, \cdots, e_t}_{\varepsilon_t})$

Let zeros$(r)$ be the position-set of the $0$'s of $r$. Similarly define ones$(r)$ and twos$(r)$.
For instance, the row $r_2$ in Table 1 has ones$(r_2) = \{5,9\}$. 
The cardinality of a row $r$, i.e. the number of characteristic vectors contained in it, clearly is

(4) \quad $|r|\quad =\quad 2^\gamma \cdot (2^{\varepsilon_1} -1) \cdots (2^{\varepsilon_t} -1)$.

Calculating the number Card$(r,k)$ of $k$-element members of $r$ is easy for say Card$(r_6, 8) =1$ and Card$(r_6, 6)=8$, but is generally harder than (4).   According to [W1, Thm.1] it holds that:

(5) \quad Calculating Card$(r,1)$, Card$(r,2), \cdots,$ Card$(r,k)$ costs $O(kw^2\log^2w)$.

We mention that merely calculating Card$(r,k)$ has the same complexity. 

Generally $\tau_k$ can be calculated as
$$\tau_k = \ds\sum_{i=1}^R \ \mbox{Card}(r_i, k),$$

where $R$ will always denote the number of (final) rows obtained by the $e$-algorithm. In our example $R = 7$ and Card$(r_i, k)$ is the entry in the $k$-th row and $i$-th column of Table 2. (Here we ignore the ``$0$-th'' row $0, \cdots, 0,1$.)

\begin{tabular}{l|c|c|c|c|c|c|c|c|c|c|c|} 
$k$& $r_1$ & $r_2$ & $r_3$ & $r_4$ & $r_5$ & $r_6$ & $r_7$ & \qquad \quad  & $\tau_k$ as sum & \qquad \quad & $f_k = {14 \choose k} - \tau_k$ \\ \hline
$0$ & 0 &0 &0 & 0 & 0 &0 &0 & &0 & & 1 \\ \hline 
$1$ & 0 &0 & 0 &0 &0 &0 &0 & & 0 & & 14 \\ \hline
$2$ & 0 &0 & 0 & 0 &0 &0 &0 && 0 & & 91 \\ \hline
$3$ & 0 &0 & 0 & 0 &0 &0 &0 && 0 & & 364\\ \hline
$4$ & 48&12&4& 0 & 2 &0 &0 & & 66 & & 935\\ \hline
$5$ &312 & 66 & 24 & 2 & 9 & 4 & 2 & & 419 & & 1583\\ \hline
$6$ & 916 & 160 & 61 & 7 & 16 & 8 & 3 & & 1171 & & 1832 \\ \hline
$7$ &1606 & 225 & 85 & 9 & 14 & 5 &1 & & 1945 & & 1487 \\ \hline
$8$ & 1868 & 202 & 70 & 5 & 6 &1 &0 && 2152 & & 851 \\ \hline
$9$ & 1509 & 119 & 34 & 1 & 1 &0 & 0 && 1664 & & 338 \\ \hline
$10$ & 858 & 45 & 9 & 0 &0 &0 &0 && 912 && 89 \\ \hline
$11$ & 339 & 10 & 1 &0 &0 &0 &0 && 350 && 14 \\ \hline
$12$ & 89 & 1 &0 &0 &0 &0 &0 && 90 & & 1 \\ \hline
$13$ & 14 &0 &0 &0 &0 &0 &0 && 14 && 0\\ \hline
$14$ & 1 & 0 &0 &0 &0 &0 &0 & & 1 & & 0 \\ \hline
& & & & & & & & & & &\\ \hline
$|r_i|$ as sum & 7560 & 840 & 288 & 24 & 48 & 18 & 6 &&  &&  \\ \hline \end{tabular}

Table 2

As expected the numbers $\tau_k$ sum up to $|{\cal S}{\cal F}'| = 8784$ (and so do the column sums $|r_i|$), and the $f_k$'s sum up to $|{\cal S}{\cal C}'| = 7600$. Of course $8784+7600 = 2^{14}$.

For each hypergraph ${\cal H}$ on $[w]$ one has ${\cal T}r({\cal H}) \neq \emptyset$ since $[w] \in {\cal T}r({\cal H})$. According to the proof of Theorem 3 of [W1] partitioning ${\cal T}r({\cal H})$ into $R> 0$ many multivalued rows costs $O(Rh^2w^2)$. 
By (5) calculating Card$(r,1), \cdots, \ \mbox{Card}(r,w)$ for any fixed row $r$ costs $O(ww^2\log^2w)$. Hence we have derived the following.

\begin{tabular}{|l|} \hline \\
{\bf Theorem 1:} With notation as above the $e$-algorithm calculates all face numbers\\
$f_k \ (1 \leq k \leq w)$ in time $O(Rh^2w^3\log^2 w)$. \\ \\ \hline \end{tabular}

If one wants $f_k$ for only {\it one} fixed $k \in [w]$ this method often works well in practise (at least for random instances) but forbids a theoretic assessment since final rows $r$ with Card$(r,k) = 0$ would be wasted\footnote{If {\it all} $f_k$ are sought, no $r$ is useless since always Card$(r,k) \neq 0$ for {\it some} $k$. See also 3.2.}  work which cannot be gauged unless such $r$ are prevented altogether.

Coming back to the $O(Rh^2w^3 \log^2w)$ bound in Theorem 1, unfortunately $R$ seems to be boundable only by $N =|{\cal S}{\cal C}|$ although in practise often $R\ll N$. Consequently $O(Nh^2w^3\log^2w)$ looks as bad as $O(2^hw)$ for inclusion-exclusion  because both $2^h< N$ and $2^h> N$ can occur. From this point of view Theorem 1 is not impressive.  However, $2^h$ is a ``real'' cost factor whereas $N$ isn't (only $R$ is). The bottom line,  one should not judge an algorithm by its worst case complexity, but by its performance in practise. 

{\bf 3.1} \ \ So let's compare the two methods (implemented with Mathematica 9.0) on some test cases. Specifically, we chose random collections of $h$ many $m$-element subsets ($=$ facets) of a $w$-set and used the $e$-algorithm ($=$ exclusion), respectively inclusion-exclusion, to calculate {\it all} face numbers $f_k \ (1 \leq k \leq w)$. The times in sec. are listed in the  last two columns of Table 3. All times below 500 sec are averages of $4$ test runs; the variance was low. Less test runs may have been averaged for larger times. For the $e$-algorithm we also list the number $R$ of final rows. It is obvious that in many cases $R$ is minuscule with respect to $N = |{\cal S}{\cal C}| > 2^m$.

Keeping $w$ fixed inclusion-exclusion scales essentially proportional to $2^h$, and so is soon outperformed as witnessed by the first three rows of Table 3. If $h$ is fixed to a small value, inclusion-exclusion scales about proportioal to $w$ (since {\it all} of $f_1, f_2, \cdots, f_w$ are required). Therefore - a small triumph - it beats exclusion which struggles to handle larger and larger multivalued rows. For instance if $h = 15$ but $w = 2000$ (and $m = 400$) then inclusion-exclusion clocks in at $1035$ sec whereas exclusion needs $15 490$ sec.

Dropping inclusion-exclusion let us focus on the comparatively benign exponential $h$-dependence of exclusion. By peculiarities of  the $e$-algorithm not further discussed, it performs the better the smaller $m/w$. Table 3 exhibits three runs of $h$ with $m/w$ fixed to $\frac{1}{3}, \frac{2}{3}, \frac{1}{6}$ respectively. In the first run (with $w = 30$) the best least square fit $y = ab^h$ to the data points $(h_1, y_1) = (15, 0.11)$ to $(h_9, y_9) = (6000, 25766)$ has $a =1566.1$ and $b =1.00048$ (as opposed to $b \approx 2$ for inclusion-exclusion). The second run (with $w = 60$) has $a = 93.3$ and $b =1.092$. The third run (with $w = 1200$) has $a = 305$ and $b = 1.144$.

\begin{tabular}{c|c|c|c|c|l}
$w$ & $h$ & $m$ & $R$ & exclusion & inclusion-exclusion \\ \hline
30 & 15 & 10 & 208 & 0.11 & 7.2 \\ \hline
30 & 16 & 10 & 218 & 0.13 & 14.5 \\ \hline
30 & 17 & 10 & 281 & 0.16 & 29.4\\ \hline
30 & 1000 & 10 & 78274 & 808 & \\ \hline 
30 & 2000 & 10 & 151193 & 3251 & \\ \hline
30 & 3000 & 10 & 218768 & 7131 & \\ \hline
30 & 4000 & 10 & 274166 & 12253 & \\ \hline
30 & 5000 & 10 & 326621 & 18683 & \\ \hline
30 & 6000 & 10 & 376290 & 25766 & \\ \hline
& & & & & \\ \hline
60 & 10 & 40 & 1750 & 2.5 & 0.4 \\ \hline 
60 & 20 & 40 & 87312 & 89 & 383 \\ \hline
60 & 30 & 40 & 929782 & 844 & \\ \hline
60 & 40 & 40 & 2302535 & 2063 & \\ \hline
60 & 50 & 40 & 10 340 983 & 8876 & \\ \hline
60 & 60 & 40 & 20 187 530 & 17956 & \\ \hline
& & & & & \\ \hline
1200 & 10 & 200 & 1825 & 47 & 5.1\\ \hline
1200 & 15 & 200 & 17245 & 395 & 163.5\\ \hline
1200 & 20 & 200 & 84018 & 1896 & 5265 \\ \hline
1200 & 25 & 200 & 271738 & 7121 & \\ \hline
1200 & 30 & 200 & 576208 & 16755 & \\ \hline
1200 & 35 & 200 & 1 145 863 & 37777 & \\ \hline
1200 & 40 & 200 & 1 931 528 & 64606 & \\ \hline
\end{tabular}

Table 3



{\bf 3.2} Rather than substituting $N$ for $R$ in Theorem 1, the experiments in 3.1 (and many others) show that the way to appreciate Theorem 1 is to view $R \ll N$ as small and to focus on how $h$ and $w$ influence the cost (namely as $h^2$ and $w^3$). 

Still, independent of significance in practice, is there a theoretic bound for the time to calculate the $f_k$'s that doesn't involve $R$ or $N$? Yes there is. Let us first consider any {\it specific} value $k \in [w]$. By (2) the cost to calculate $f_k$ is bound by the complexity to calculate $\tau_k$, which is known to be fixed-parameter-tractable in time $O(d^{2k+1} h)$. See [FG, p.358] where the fixed parameters are $k$ and $d: = \max \{|H_i|: \ 1 \leq i \leq h\}$ with $H_i = W \backslash F_i$. Observe that using this method to calculate {\it all} $f_k$'s costs $O(d^{2w+1} hw)$ which wouldn't qualify as fixed-parameter-tractable because $w$ (which cannot be ``fixed'') appears in the exponent.  It would be interesting to pit this method against the one from Theorem 1.

\section{Partitioning a simplicial complex, given its facets}

The {\it face lattice} of a simplicial complex ${\cal S}{\cal C}$ based on $W$ is ${\cal L}: = {\cal S}{\cal C} \cup \{W\}$, ordered by inclusion. For all $X, Y \in {\cal L}$ the meet $X \wedge Y$ is $X \cap Y$; the join $X \vee Y$ is $X \cup Y$ provided $X \cup Y\in {\cal S}{\cal C}$, and $W$ otherwise. According to [KP] the {\it combinatorial face lattice enumeration problem} is the following: Given the facets $F_1, \ldots, F_h$, calculate the diagram of ${\cal L}$, which thus entails a listing of all covering pairs of faces. Letting $N = |{\cal L}|$ (or $N = |{\cal S}{\cal C}| = |{\cal L}|-1$) this can be achieved [KP, Thm.5] in time $O(N\min (h,w)(|F_1| + \cdots +|F_h|)$ which is $O(Nh^2w)$ or $O(Nhw^2)$ depending on whether $h \geq w$ or $h \leq w$. The space required is $O(N\min (h,w))$.

What about merely enumerating ${\cal L}$ and dispensing with its edges? The naive way\footnote{By this we mean enumerating in turn ${\cal P}(F_1)$ to ${\cal P}(F_h)$ and checking for each $X \in {\cal P}(F_i)$ whether it is among the faces retrieved so far.} takes time $O(N^2hw)$. But ${\cal L}$ being a closure system one can do it in time $O(Nhw^2)$ and space $O(hw)$. Namely, let ${\cal L} \subseteq {\cal P}([w])$ be {\it any} closure system which is given by its associated closure operator $c: \ {\cal P}([w]) \ra {\cal P}([w])$. If computing an arbitrary closure $c(X)$ costs at most $f(w)$   then enumerating ${\cal L}$ can be done\footnote{This follows at once from [GR] and earlier work of Ganter, albeit the $O(Nwf(w))$ bound is not explicitly stated in [GR]. More specifically, the elements of ${\cal L}$ are output one by one with {\it delay} $O(wf(w))$.} in time $O(Nwf(w))$. It is clear that in our situation $f(w) = O(hw)$.

In this Section we present a method, very much different from [KP] or [GR], which delivers ${\cal S}{\cal C}$ in a compact way, i.e. not one by one.

We henceforth use the symbol $\uplus$ for {\it disjoint} union. For $1 \leq p \leq h$ let ${\cal S}{\cal C}_p \subseteq {\cal S}{\cal C}$ be the simplicial complex generated by $F_1, \cdots, F_p$, so ${\cal S}{\cal C}_p = {\cal P}(F_1) \cup \cdots \cup {\cal P}(F_p)$. By induction assume that ${\cal S}{\cal C}_p = r_1 \uplus r_2 \uplus \cdots \uplus r_m$ with $\{0,1,2,e\}$-valued rows $r_i$. The basic procedure to extend this to a representation of ${\cal S}{\cal C}_{p+1}$ is as follows. 
We iteratively shrink
${\cal F}_{p+1} [0] = {\cal P}(F_{p+1})$ (viewed as $\{0,2\}$-valued row) in order to make it disjoint from ${\cal S}{\cal C}_p$.  Thus ${\cal F}_{p+1}[0] \supseteq {\cal F}_{p+1} [1] \supseteq {\cal F}_{p+1}[2] \supseteq \cdots$ until ${\cal F}_{p+1} [p]$ is such that

(6) \qquad ${\cal S}{\cal C}_{p+1} = {\cal S}{\cal C}_p \cup {\cal F}_{p+1} [0] = {\cal S}{\cal C}_p \uplus {\cal F}_{p+1} [p]$.

Since ${\cal F}_{p+1}[p]$  itself arises as disjoint union of multivalued rows, the induction hypothesis will carry over from ${\cal S}{\cal C}_p$ to ${\cal S}{\cal C}_{p+1}$.
It is clear that ${\cal F}_{p+1}[p]$ satisfies (6) if we set

(7) \qquad ${\cal F}_{p+1} [i] : = \{X \in {\cal F}_{p+1} [i-1]: \ X \not\in {\cal P}(F_i) \} \quad (1\leq i \leq p)$.

For the simplicial complex ${\cal S}{\cal C} = {\cal S}{\cal C}_6$ from Section 2 the procedure unfolds as follows. (See Table 4; its rows $r_1$ to $r_7$ have nothing to do with $r_1$ to $r_7$ in Table 1.)

\begin{tabular}{l|c|c|c|c|c|c|c|c|c|c|c|c|c|c|l}
 & 1 & 2 & 3& 4 & 5 &6 & 7 & 8 & 9 & 10 & 11 & 12 & 13 & 14 & \\ \hline
 $r_1=$ & 2 & 2 & 0 & 0 & 2 & 2 & 2 & 2 &  0 & 2 & 2 & 2 & 2 & 2 & ${\cal F}_1[0]$\\ \hline
 & 2 & 2 & 2 & 2 & 0 & 2 & 2 & 2 & 2 & 0 & 2 & 2 & 2 & 2 & ${\cal F}_2[0]$ \\ \hline
 $r_2=$ & 2 & 2 & $e_1$ & $e_1$ & 0 & 2 & 2 & 2 & $e_1$ & 0 & 2 & 2 & 2 & 2 & ${\cal F}_2[1]$\\ \hline
 & 2 & 2 & 2 & 2 & 2 & 0 & 0 & 2 &  2 & 2 & 0 & 0 & 2 & 2 & ${\cal F}_3[0]$\\ \hline 
 $r_3=$ & 2 & 2 & $e_1$ & $e_1$ & $e_2$ & 0 & 0 & 2 & $e_1$ & $e_2$ & 0 & 0 &  2 & 2 & ${\cal F}_3[2]$ \\ \hline
 & 2 & 2 & 2 & 2 &  2 & 2& 2 & 0 & 2 & 2 & 2 & 2 & 0 & 0 & ${\cal F}_4[0]$ \\ \hline
 $r_4=$ & 2 & 2 & $e_1$ & $e_1$ & $e_2$ & $e_3$ & $e_3$ & 0 & $e_1$ & $e_2$ & $e_3$ & $e_3$ & 0 & 0 & ${\cal F}_4[3]$\\ \hline
 & 0  & 0 & 0 & 0 & 0 & 0 & 0 & 0 & 2 & 2 & 2 & 2 & 2 & 2 & ${\cal F}_5[0]$ \\ \hline
 $r_5=$ & 0 & 0 & 0 & 0 & 0 & 0 & 0 & 0 & 1 & 1 & $e_3$ & $e_3$ & $e_4$ & $e_4$ & ${\cal F}_5[4]$\\ \hline
 & 2 & 2 & 0 & 0 & 0 & 2 & 2 & 0 & 2 & 2 & 2 & 0 &0 & 2 & ${\cal F}_6[0]$\\ \hline
 & 2 & 2 & 0 & 0 & 0 & $e_3$ & $e_3$ & 0 & 1 & 1 & $e_3$ &  0 & 0 & 1 & ${\cal F}_6[4]$\\ \hline
 $r_6=$ & 2 & 2 & 0 & 0 & 0 & $e_3$ & $e_3$  & 0 & 1 & 1 & 2 & 0 & 0 & 1 & \\ \hline
 $r_7=$ & $e_5$ & $e_5$ & 0 & 0 & 0 & 0 & 0 & 0 & 1 & 1 & 1 & 0 & 0 & 1 & \\ \hline
 \end{tabular}
 
 Table 4

Starting with
$$r_1: = {\cal F}_1[0] = {\cal P}(F_1) = (2,2,0,0,2,2,2,2,0,2,2,2,2,2)$$
the only way for $X \in {\cal F}_2[0] = {\cal P}(F_2)$ not to be a member of $r_1$ is to have $X \cap \ \mbox{zeros}(r_1)  = X \cap \{3,4,9\} \neq \emptyset$. Hence
$$r_2: = {\cal F}_2[1] = (2,2,e_1, e_1, 0,2,2,2,e_1, 0,2,2,2,2).$$
Similarly ${\cal F}_3[1]$ arises from ${\cal F}_3[0]$ by putting $e_1 e_1 e_1$ on positions $3,4,9$; and $r_3: = {\cal F}_3[2]$ arises from ${\cal F}_3[1]$ by putting $e_2e_2$ on
zeros$({\cal P}(F_2)) = \{5,10\}$. So far
$${\cal S}{\cal C}_3 = {\cal P}(F_1) \cup {\cal P}(F_2) \cup {\cal P}(F_3) = r_1 \uplus r_2 \uplus r_3,$$
and we continue the same way up to ${\cal F}_6[0]$ for which the subset ${\cal F}_6[5]$ can no longer be represented as a {\it single} $\{0,1,2,e\}$-valued row. Instead we note that the partition
$${\cal F}_6[4]\quad =\quad \{X \in {\cal F}_6[4]: X \cap\{6,7\} \neq \emptyset \} \quad \uplus \quad \{X \in {\cal F}_6[4]: \ X \cap \{6,7\} = \emptyset \}$$
displays as follows in terms of multivalued rows:

\begin{tabular}{l|c|c|c|c|c|c|c|c|c|c|c|c|c|c|c|} 
& 1 & 2 & 3 & 4 & 5 & 6 & 7 & 8 & 9 & 10 & 11 & 12 & 13 & 14 \\ \hline
${\cal F}_6[4] =$ & 2 & 2 & 0 & 0 & 0 & $e_3$ & $e_3$ & 0 & 1 & 1 & $e_3$ & 0 & 0 & 1 \\ \hline
$r^+ =$ & 2 & 2 & 0 & 0 & 0 & ${\bf e}_3$ & ${\bf e}_3$ & 0 & 1 & 1  & ${\bf 2}$ & 0 & 0 & 1 \\ \hline
$r^-  =$ & 2 & 2 & 0 & 0 & 0 & ${\bf 0}$ & ${\bf 0}$ & 0 & 1 & 1  & ${\bf 1}$ & 0 & 0 & 1 \\ \hline \end{tabular}

The rows $r^-$ and $r^+$ are the {\it candidate sons} of $r$. Generally  [W2, Sec. 4], with $t$ as in (1), there can be up to $t$ candidate sons and they are always such that the pending constraint can be smoothly imposed upon them. Here $r_6 : = r^+ \subseteq {\cal F}_6[5]$ and $r^- \cap {\cal F}_6[5]$ can be written as a multivalued row (namely $r_7$ in Table 4). Thus ${\cal F}_6[5] = r_6 \uplus r_7$, and so ${\cal S}{\cal C}_6 = r_1 \uplus \cdots \uplus r_7$. The latter equality is also supported by the calculation
$$|r_1| + \cdots + |r_7|\quad = \quad 2048 + 3584 + 672+ 1260 + 9 + 24 + 3 \quad =\quad 7600$$
where 7600 is the number obtained in Section 2.

 In general ${\cal F}_p [i-1]$ may be $\rho_1 \uplus \rho_2 \uplus \cdots$ and upon enforcing each row's disjointness from ${\cal P}(F_i)$ each $\rho_j$ may in turn decay into a disjoint union of rows. The union of the latter rows yields ${\cal F}_p[i]$. 
 
 \begin{tabular}{|l|} \hline \\
 {\bf Theorem 2:} Given are the facets $F_1, \cdots, F_h \subseteq [w]$ of a simplicial
 complex ${\cal S}{\cal C}$. Mentioned \\
  method (based on the $e$-algorithm) to partition ${\cal S}{\cal C}$ into $R$ pieces
 costs $O(Rh^3w^2)$.\\ \\ \hline \end{tabular}

{\it Proof:} Since $\{F_1, \cdots, F_h\}$ is an antichain (i.e. $F_i \not\subseteq F_j$ for $i \neq j$), no ${\cal P}(F_{p+1}) = {\cal F}_{p+1}[0]$ is contained in ${\cal S}{\cal C}_p= {\cal P}(F_1) \cup \cdots \cup {\cal P}(F_p)$. Hence ${\cal F}_{p+1}[p] \neq \emptyset$ in (6). Representing ${\cal F}_{p+1}[p]$ by a disjoint union of $R_{p+1}> 0$ many $\{0,1,2,e\}$-valued rows entails imposing on ${\cal F}_{p+1}[0]$ the $p$ constraints $X \cap H_i \neq \emptyset$ (where $H_i: = W \setminus F_i$ for $1 \leq i \leq p$), and this costs $O(R_{p+1} p^2w^2)$ before Theorem 1. In view of $R=1+R_2 + \cdots R_h$ the claim now follows from
$$\ds\sum_{p=1}^{h-1} O(R_{p+1} p^2w^2) = \ds\sum_{p=1}^{h-1} O(Rp^2w^2) = O(Rh^3w^2).$$
\hfill $\square$

Compare $O(Rh^3w^2)$ in Theorem 2 with the $O(Nhw^2)$ algorithm of [GR] for outputting the faces {\it one by one}. Albeit $hw^2 < h^3w^2$, this is outweighed by  $N \gg R$. Also compare the $O(Rh^3w^2)$ cost of {\it partitioning} with the $O(Rh^2w^3\log^2w)$ cost for merely calculating the {\it face numbers} (Theorem 1).

{\bf 4.1} \ \ One use of partitioning a simplicial complex ${\cal S}{\cal C}$ is optimization. Given a target function $f: \ [w] \ra \Z$, extend it to $f : \ {\cal S}{\cal C} \ra \Z$ as usual by putting $f(X):= \sum \{f(i): \ i \in X\}$. It is easy to find
$$\mu ({\cal S}{\cal C}): = \max \{f(X): X \in {\cal S}{\cal C} \}$$
when ${\cal S}{\cal C} = r_1 \uplus \cdots \uplus r_R$ is displayed as a disjoint union of multivalued rows. Namely, 
$$\mu ({\cal S}{\cal C}) = \max \{\mu (r_i): 1 \leq i \leq R\}$$
 where $\mu (r_i)$ is readily obtained, as can be gleaned from this example with $w = 15$:

\begin{tabular}{l|c|c||c|c||c|c|c|c||c|c|c|c||c|c|c|} \hline
$r=$ &  0 & 0 & 1 & 1 & 2 & 2 & 2 & 2 & $e_1$ & $e_1$ & $e_1$ & $e_1$ & $e_2$ & $e_2$ & $e_2$ \\ \hline
$f=$ & 8 & 5 & ${\bf -6}$ & ${\bf -2}$ & ${\bf 3}$ & $-5$ & $-7$ & ${\bf 8}$ & $-4$ & $-2$ & ${\bf 5}$ & ${\bf 4}$ & $-10$ & ${\bf -8}$ & $-9$ \\ \hline \end{tabular}

$\mu (r) = (-6) + (-2) + 3 + 8+5 + 4 + (-8)  =4$. More complicated target functions, such as quadratic forms, are also amenable to this approach.

{\bf 4.2} Similar to 3.2 one may wonder: What if we wish to enumerate (in ordinary set notation) all $f_k$ faces of {\it specific} cardinality $k$? Let us generalize this to the task of enumerating the set $I_k (P, \leq)$ of all $k$-element order ideals of a $w$-element poset $(P, \leq)$; it is unlikely that the cost for the special case of simplicial complexes is much lower. If $|P| < k$ then $I_k (P, \leq) = \emptyset$. Otherwise pick any $u \in P$ and split $I_k (P, \leq)$ into the subsets
$$S^+ = \{X \in I_k (P, \leq) : u \in X\} \quad \mbox{and} \quad S^- = \{X \in I_k(P, \leq): \ u \not\in X\}.$$
At least one of $S_+$ and $S_-$ is nonempty and there are obvious posets $(P^+, \leq)$ and $(P^-, \leq)$ such that $S^+$ and $S^-$ bijectively correspond to $I_k (P^+, \leq)$ and $I_k(P^-, \leq)$ respectively. Iterate this for $(P^+, \leq)$ and $(P^-, \leq)$ and so forth. This gives rise to a binary recursion tree with root $(P, \leq)$ whose $f_k$ leaves are the $k$-element order ideals of $P$ (we use the same notation $f_k$ as for simplicial complexes). It is easy to see that all of this costs $O(f_kw^3)$. We mention that the whole procedure requires to carry along the sets $y\uparrow = \{v \in P: v \geq y\}$ and $y \downarrow = \{v \in P: v \leq y\}$ for all $y \in P$. Updating these sets when switching from $P$ to $P^- = P \backslash (u \uparrow)$ costs $O(w^2)$. 

A quite different $O(f_k w^3)$ algorithm to enumerate the $k$-element order ideals is given in [W4]. Its advantage,  as in the present article (but tackling a specific $k$), is that the sought objects are output in clusters, not one by one.

\section{Partitioning a simplicial complex, given its minimal nonfaces}

Suppose ${\cal S}{\cal C}$ is a simplicial complex whose minimal non-faces $G_1, \cdots, G_q$ are known, i.e. the $G_i$'s are the generators of the set filter ${\cal S}{\cal F} = {\cal P}(W) \setminus {\cal S}{\cal C}$. Then ${\cal S}{\cal C}$ consists of all {\it noncovers} $X$ of $\{G_1, \cdots, G_q\}$ i.e. of all sets $X \subseteq W$ such that

(8) \qquad $(\forall 1 \leq i \leq q) \quad X \not\supseteq G_i$.

Because (8) is equivalent to

(9) \qquad $(\forall 1 \leq i \leq q) \quad (W \setminus X) \cap G_i \neq \emptyset$,

one can employ\footnote{In Section 3 we generated ${\cal S}{\cal F}$ (which suffices for the {\it face numbers} of ${\cal S}{\cal C}$) by imposing $h$ constraints $(W \setminus F_i) \cap X \neq \emptyset$. Here in contrast we generate ${\cal S}{\cal C}$ by imposing $q$ constraints $(W 
\setminus X ) \cap G_i \neq \emptyset$. Both tasks can be achieved by the $e$-algorithm. However, as mentioned in previous publications, for tasks such as the latter it is more succinct to use a dual version of the $e$-algorithm, called $n$-{\it algorithm}. Its output is a disjoint union of $\{0,1,2,n\}$-valued rows where the wildcard $nn \cdots n$ means {\it at least one $0$ here}. The statement of Theorem 3 is mentioned in [W2, p.64], where the ``negative clause base'' is just our set of minimal nonfaces.} the transversal $e$-algorithm to enumerate all noncovers of $\{G_1, \cdots, G_q\}$, and whence enumerate ${\cal S}{\cal C}$; this costs $O(Rq^2w^2)$ where $R$ is the number of final rows.
Let us summarize.

\begin{tabular}{|l|} \hline \\
{\bf Theorem 3:} Given are the minimal non-faces $G_1 \cdots, G_q \subseteq [w]$ of a simplicial
complex\\
 ${\cal S}{\cal C}$. 
Mentioned method (based on the $n$-algorithm) to partition ${\cal S}{\cal C}$ 
into $R$ pieces costs\\
 $O(Rq^2w^2)$. \\ \\ \hline \end{tabular}
 
 For example, for the simplicial complex ${\cal S}{\cal C}'$ from Section 2 one computes 74 minimal non-faces $G_i$. Applying the noncover $n$-algorithm to the $G_i$'s displays ${\cal S}{\cal C}'$ as a disjoint union of 37 $\{ 0,1,2,n\}$-valued rows. 
As another example, let ${\cal S}{\cal C}$ be the simplicial complex of all independent sets of a matroid ${\bf M}$ on $[w]$, given by its minimal non-faces (the {\it circuits} of ${\bf M}$). Then ${\cal S}{\cal C}$ can be partitioned with the noncover $n$-algorithm. In this particular case the facets of ${\cal S}{\cal C}$ (the {\it bases} of ${\bf M}$) are equicardinal. In a similiar vein, applying the $n$-algorithm to the so called {\it broken} circuits of ${\bf M}$ yields a simplicial complex whose face numbers coincide (up to sign) with the coefficients of the chromatic polynomial of ${\bf M}$. If all $G_i$'s are of cardinality 2 then ${\cal S}{\cal C}$ is the set of all anticliques of some obvious graph. In this case the $n$-algorithm can be fine-tuned. Some of these and other topics will be expanded elsewhere. One of them will be previewed in a bit more detail in Section 6.

Unfortunately, if the minimal non-faces $G_i$ are unknown, it it tough finding them. Specifically, calculating the $G_i$'s from the facets $F_i$ amounts to the hard problem [EMG] of dualizing a monotone Boolean function. 
It may thus be faster to confront the $\frac{1}{2}(h^2-h)$ many constraints from Section 4 and get ${\cal S}{\cal C}$ as disjoint union of $\{0,1,2,e\}$-valued rows rather than the $\{0,1,2,n\}$-valued rows of Theorem 3. Either $h$ or $q$ can be much larger than the other and usually one cannot tell in advance.

\section{Application to inclusion-exclusion}

As another target of the $n$-algorithm let's see how the calculations for inclusion-exclusion can be sped up. Say $a_i \ (1 \leq i \leq m)$ are any $m$ properties potentially applying to anyone of $N_0$ objects. If $N$ is the number of objects enjoying all properties and if say $N (\ol{a}_2 \ol{a}_4)$ is the number of objects violating $a_2$ and $a_4$, then for say $m = 4$:

(10) \quad $N = N_0 - \sum N(\ol{a}_i) + \sum N (\ol{a}_i \ol{a}_j) - \sum N (\ol{a}_i \ol{a}_j \ol{a}_k) + N (\ol{a}_1 \ol{a}_2 \ol{a}_3 \ol{a}_4)$.

Often many terms $N (\ol{a}_i\ol{a}_j \cdots)$ are zero, and obviously the corresponding sets $\{ \ol{a}_i, \ol{a}_j, \cdots \}$ constitute a set filter ${\cal S}{\cal F} \subseteq {\cal P}(\{\ol{a}_1, \cdots \ol{a}_m \} )$, which may be coined the {\it irelevant set filter}. Suppose again $m = 4$ and ${\cal S}{\cal F}$ is generated by $\{1, 4\} (: = \{ \ol{a}_1, \ol{a}_4\}), \{3, 4\}, \{1, 2, 3\}$. The $n$-algorithm delivers the complimentary simplicial complex ${\cal S}{\cal C}$ as a disjoint union of $\{0, 1, 2, n\}$-valued rows, in our case $r_1 = (n, n, 1, 0)$ and $r_2 = (n, 2, 0, n)$. Thus $r_1$ provides (if we let in turn $nn = 00, 01, 10)$ the nonzero terms $N(\ol{a}_3), N(\ol{a}_2 \ol{a}_3), N(\ol{a}_1 \ol{a}_3)$. Similarly $r_2$ provides the remaining six nonzero terms $N(\cdots)$ in the expansion (10). It is interesting to pit this method against {\tt SatisfiabilityCount} when the task is to count the models of a Boolean function given in disjunctive normal form [W3, Section 6].

The approach to pack ${\cal S}{\cal C}$ into multivalued rows is particularly appealing if each nonzero $N(\ol{a}_i \ol{a}_j \cdots )$ is merely a function $g(k)$ of the {\it cardinality} $k = |\{ \ol{a}_i, \ol{a}_j, \cdots \} |$. That sometimes happens when the properties $a_i$ are of a symmetric kind. Letting $f_k$ be the number of $k$-faces of ${\cal S}{\cal C}$ (computed as in Theorem 2), formula (10) then improves to

(11) \qquad $N = N_0 + \ds\sum_{k=1}^m (-1)^k g(k) f_k$.

Formula (11) e.g. applies to the count of constrained permutations [W3, Section 3], respectively integer partitions [W3, Section 5]. In both cases the effort to find the minimal non-faces of ${\cal S}{\cal C}$ (an issue raised at the end of Section 5) is moderate.

Sections 7 and 8 concern applications of the scenario where all facets are given (as opposed to the  minimal non-faces).

\section{Application to combinatorial commutative algebra}

Following [MS] in Section 7 we adopt the notation $\Delta$ instead of ${\cal S}{\cal C}$.
We shall touch upon the $h$-polynomial (7.1), reduced homology groups (7.2), and the link of a face (7.3). The deliberations in 7.1 and 7.2 are less ``ad hoc'' than the ones in 7.3. Still, 7.3 gives some additional flavor of how the use of wildcards may benefit computational commutative algebra. 

{\bf 7.1} As to the $h$-polynomial 
$$h(t) : = \ds\sum_{i=1}^w f_{i-1} t^i (1-t)^{w-i} = h_0 + h_1t + \cdots + h_wt^w$$
of $\Delta$, it is immediately computed  from the face numbers $f_i$ of $\Delta$ and has many applications [S]. Here {\it all} $f_i$ are required and they\footnote{In combinatorial commutative algebra an $i$-face is defined as having cardinality $i+1$ (but {\it dimension} $i$). We stick to our definition of cardinality $i$.} can be calculated from the facets $F_i \in \Delta$ as in Section 3. 

{\bf 7.2}  The $i$-th {\it reduced homology} $\tilde{H}_i (\Delta; K)$ over the field $K$ is defined as follows. Let $V_j$ be a vector space which has all $j$-faces as a $K$-basis. For certain $K$-linear {\it boundary maps} $\delta_i$ and $\delta_{i+1}$ of type
$$V_{i+1} \quad {\delta_{i+1} \atop \longrightarrow} \quad V_i \quad {\delta_i \atop \longrightarrow} \quad V_{i-1}$$
one can show that im$(\delta_{i+1}) \subseteq \ \mbox{ker} (\delta_i)$. Hence one can define the quotient $K$-vector space

(12) \qquad $\tilde{H}_i (\Delta; K) \ : = \ \mbox{ker}(\delta_i) / \mbox{im}(\delta_{i+1})$.

Often only the {\it dimension} of $\tilde{H}_i(\Delta; K)$ is at stake. This being
$$\begin{array}{lll}
 & & \mbox{dim}(\mbox{ker}(\delta_i)) - \ \mbox{dim}(\mbox{im} (\delta_{i+1})) \\
 \\
 & =& \mbox{dim}(V_i) - \ \mbox{dim}(\mbox{im} (\delta_i)) - \ \mbox{dim} (\mbox{im}(\delta_{i+1})), \end{array}$$
 it suffices to show how to calculate $\mbox{dim}(\mbox{im}(\delta_i))$. The latter is the rank of a $w$ by $f_i$ matrix $M$ whose columns are the ``signed'' supports (i.e. having alternating entries $\pm 1$) of the $i$-faces; see e.g. [MS, Example 1.18]. The $i$-faces can be delivered compactly encoded as shown in Section 4. Calculating rank$(M)$ is subtle and depends a lot on $K$, but for some researchers it may boil down to a hardwired command in their favorite programming language.
 
 {\bf 7.3} According to [MS, p.17] the {\it link} of a face $X$ of a simplicial complex $\Delta$ is
 
 (13) \qquad link$_{\Delta}(X) : = \{Y \in \Delta: \ Y \cup X \in \Delta \ \mbox{and} \ Y \cap X = \emptyset \}$, 
 
 i.e. the set of faces that are disjoint from $X$ but whose union with $X$ is small enough to stay in $\Delta$. Of course $\Delta (X):=$ link$_{\Delta}(X)$ is a simplicial complex itself. Simplicial complexes of type link$_{\Delta}(X)$ occur in many situations, e.g. when defining {\it Betti} {\it numbers}, or in Reisner's criterion for a Stanley-Reisner ring to be Cohen-Macaulay.
 
 We close this Section by indicating how link$_\Delta(X)$ can be partitioned into multivalued rows. The quick answer is that the facets of link$_\Delta(X)$ are the maximal members $G_j$ of the set family 
 $$\{F_i \setminus X: \ 1 \leq i \leq h, \ X \subseteq F_i\}$$
 and so applying the method of Section 4 to the sets $G_j$ does the job. However,  an existing partition of $\Delta$ can be exploited to get a partition of link$_\Delta(X)$ faster. To do so define
 $$\begin{array}{lll}\mbox{Disjoint}(\Delta, X) & : = & \{Y \in \Delta: \ Y \cap X = \emptyset \}, \\
 \mbox{Minus}(\Delta, X) & : = & \{Z \setminus X: \ Z \in \Delta, \ X \subseteq Z\}. \end{array}$$
 Albeit a member of Minus$(\Delta, X)$ needs not be a face of ${\cal S}{\cal C}$, it is clear that
 
 (14) \quad link$_\Delta(X) = \ \mbox{Disjoint}(\Delta, X) \cap \ \mbox{Minus}(\Delta, X)$.
 
 As will be seen, $\Delta = r_1 \uplus \cdots \uplus r_R$ readily spawns partitions
 $$\mbox{Disjoint}(\Delta, X) \quad = \quad \rho_1 \uplus \cdots \uplus \rho_\alpha,  \quad \mbox{and} \quad  \mbox{Minus}(\Delta, X) \quad = \quad \sigma_1 \uplus \cdots \uplus \sigma_\beta,$$
 and so it follows from (13) that
 
 (15) \quad link$_\Delta (X) = \uplus \{\rho_i \cap \sigma_j : \ i \in [\alpha], \ j \in [\beta] \}$.
 
 It remains to write all nonempty $\rho_i \cap \sigma_j$ as a disjoint union of multivalued rows.
 
 To fix ideas, take the simplicial complex $\Delta = {\cal S}{\cal C}'$ from Section 2 which is $\Delta = r_1 \uplus \cdots \uplus r_7$ according to Table 4. For the face $X = \{6,7,10,11\}$ one verifies that Disjoint$(\Delta, X) = \ol{r}_1 \uplus \ol{r}_2 \uplus \ol{r}_3 \uplus \ol{r}_4$ and Minus$(\Delta, X) = r'_1 \uplus r'_4 \uplus r'_6$ (see Table 5).
  Seven of the $4\cdot 3=12$ intersections $\ol{r}_i \cap r'_j$ are empty, the other five happen to be expressible as {\it single} multivalued rows as shown in Table 5. According to (15), link$_\Delta(X)$ is the disjoint union of these rows.
 
 In general the intersection of two multivalued row $\rho_i$ and $\sigma_j$ is handled by imposing the $e$-bubbles (and $0$'s and  1's) of one row upon the other. If merely $|\rho_i \cap \sigma_j|$ is required this can be obtained with inclusion-exclusion. Specifically, let the row with the fewer $e$-bubbles have $m$ $e$-bubbles. Then $m$ is the number of properties $a_i$ to be dealt with (see also Sec.6).
 
 \begin{tabular}{l|c|c|c|c|c|c|c|c|c|c|c|c|c|c|} 
  & 1 & 2& 3& 4 &5 & 6 & 7 &8 & 9 & 10 & 11& 12 & 13 & 14\\ \hline
  $\ol{r}_1=$ & 2 & 2 & 0 & 0 & 2 & ${\bf 0}$ & ${\bf 0}$ & 2 & 0 & ${\bf 0}$ & ${\bf 0}$ & 2 & 2 & 2 \\ \hline
  $\ol{r}_2=$ & 2 & 2 & $e_1$ & $e_1$ & 0 & ${\bf 0}$ & ${\bf 0}$ & 2 & $e_1$ & ${\bf 0}$ & ${\bf 0}$ & 2 & 2 & 2\\ \hline
  $\ol{r}_3=$ & 2 & 2 & $e_1$ & $e_1$ & 1 & ${\bf 0}$ & ${\bf 0}$ & 2 & $e_1$ & ${\bf 0}$ & ${\bf 0}$ & 0 & 2 & 2  \\ \hline
  $\ol{r}_4=$ & 2 & 2 & $e_1$ & $e_1$ & 1 & ${\bf 0}$ & ${\bf 0}$ & 0 & $e_1$ & ${\bf 0}$ & ${\bf 0}$ & 1 & 0 & 0 \\ \hline
  & & & & & & & & & & & & & & \\ \hline
  $r'_1=$ & 2 & 2 & 0 & 0 & 2 & ${\bf 0}$ & ${\bf 0}$ & 2 & 0 & ${\bf 0}$ & ${\bf 0}$ & 2 & 2 & 2 \\ \hline
  $r'_4=$ & 2 & 2 & $e_1$ & $e_1$ & 2 & ${\bf 0}$ & ${\bf 0}$ & 0 & $e_1$ & ${\bf 0}$ & ${\bf 0}$ & 2 & 0 & 0\\ \hline
  $r'_6=$ & 2 & 2 & 0 & 0 & 0 &  ${\bf 0}$ & ${\bf 0}$ &  0 & 1 & ${\bf 0}$ & ${\bf 0}$ & 0 & 0 & 1\\ \hline
   & & & & & & & & & & & & & & \\ \hline
   $\ol{r}_1 \cap r'_1=$ & 2 & 2 & 0 & 0 & 2 & 0 & 0 & 2 & 0 & 0 & 0 & 2 & 2 & 2\\ \hline
   $\ol{r}_2 \cap r'_4 =$ & 2 & 2 & $e_1$ & $e_1$ & 0 & 0 & 0 &0 & $e_1$ & 0 & 0 & 2 &0 & 0 \\ \hline
   $\ol{r}_2 \cap r'_6=$ & 2 & 2 & 0 & 0 & 0 & 0 & 0 & 0 & 1 & 0 & 0 & 0 & 0 & 1\\ \hline
   $\ol{r}_3 \cap r'_4 =$ & 2 & 2 & $e_1$ & $e_1$ & 1 & 0 & 0 & 0 & $e_1$ & 0 & 0 & 0 & 0 & 0 \\ \hline
   $\ol{r}_4 \cap r'_4=$ & 2 & 2 & $e_1$ & $e_1$ & 1 & 0 & 0 & 0 & $e_1$ & 0 & 0 & 1 & 0  &0 \\ \hline
      \end{tabular}
      
      Table 5

\section{Application to Frequent Set Mining}

A property $\pi$ defined on the subsets $X \subseteq W$, is called {\it monotone} if with $X$ each subset of $X$ enjoys $\pi$. Evidently the set of all $X$'s that enjoy $\pi$ constitutes a simplicial complex.
As to one particular monotone property to be dealt with, fix ${\cal D} \subseteq {\cal P}(W)$ and an integer $s \geq 1$. Then $X$ is called $s^+$-{\it frequent} (with respect to ${\cal D}$) if $X$ is a subset of at least $s$ members $T \in {\cal D}$. For instance, let $W = [4]$ and ${\cal D}:= \{\{1,2,3\}, \{1,2,4\}, \{3,4\}\}$. Then the simplicial complex ${\cal S}{\cal C}$ of all $2^+$-frequent sets is
${\cal S}{\cal C} = \{ \emptyset, \{1\}, \{2\}, \{3\}, \{4\}, \{1,2\}\}.$

In the framework of {\it Frequent Set Mining} (FSM) one refers to ${\cal D}$ as the {\it database} and to its elements $T \in {\cal D}$ as {\it transactions}. For instance, the transactions could comprise the items bought by  customers in a supermarket during a specific period of time. Hence an {\it itemset} $X$ is $s^+$-frequent if its items have been bought {\it together} at least $s$ times. We write ${\cal F}r(s^+)$ for the simplicial complex of all $s^+$-frequent sets. If ${\cal F}r(s^+)$ is small, its faces can be enumerated one by one, and in the early days of FSM the so called {\it A priori} algorithm got famous for doing just that: it even made it to the Top Ten Algorithms in Data Mining [WK, ch.4]. 

For large ${\cal F}r(s^+)$ it may still be desirable to have all of ${\cal F}r(s^+)$ available, but then ${\cal F}r(s^+)$ must be encoded in compact form, somehow. As seen in Sections 4 and 5, this can be done if the facets of ${\cal F}r(s^+)$ are available\footnote{We mention in passing that the facets can be calculated from ${\cal D}$ with another standard FSM algorithm, called {\it Dualize and Advance}.}. In 8.1 we discuss how some telling probabilities can be computed from the facets of ${\cal F}r(s^+)$. Subsection 8.2 concerns so called {\it closed} frequent sets.

{\bf 8.1} \quad We shall focus on the calculation of certain numbers attached to ${\cal F}r(1^+), {\cal F}r(2^+)$, and so on. Put ${\cal D} = \{T_1, T_2, \cdots, T_m \}$. For simplicity assume that ${\cal D}$ is an antichain. Then the facets of ${\cal S} {\cal C} : = {\cal F}r(1^+)$ are just $T_1$ up to $T_m$. We are going to slice ${\cal S} {\cal C}$ in two ways. First, cardinality-wise as

(16) \quad ${\cal S} {\cal C} \quad = \quad {\cal S}{\cal C}[0] \uplus {\cal S}{\cal C}[1] \uplus \cdots \uplus {\cal S}{\cal C}[\gamma]$,

where $\gamma$ is the maximum cardinality of a facet. Second, frequency-wise as

(17) \quad ${\cal S}{\cal C} \quad = \quad {\cal F}r(1) \uplus {\cal F}r(2) \uplus \cdots \uplus {\cal F}r(m)$,

where ${\cal F}r(s) : = {\cal F}r(s^+) \setminus {\cal F}r((s+1)^+)$ is the family of all $s$-{\it frequent}\footnote{Be aware that in the FSM literature usually ``$s$-frequent'' corresponds to our $s^+$-frequent.} sets in the sense that they occur in {\it exactly} $s$ transactions. Furthermore, let
$$fr (s,k) \quad : = \quad |{\cal F}r(s) \cap {\cal S}{\cal C}[k]| \quad (0 \leq k \leq \gamma, \quad 1 \leq s \leq m)$$
be the number of $k$-element subsets of $W$ which are $s$-frequent. These numbers can be computed as

\hspace*{2.6cm} $fr(s,k) \quad =\quad fr(s^+, k) - fr ((s+1)^+, k)$

where 
$$fr(s^+,k) \quad := \quad \mbox{number of} \ s^+\mbox{-frequent sets of cardinality} \ k.$$
Dually to $f(s^+, k)$ we define 
$$f(s, k^+) \quad :=  \quad \mbox{number of} \  s\mbox{-frequent sets of cardinality} \  \geq k.$$
Obviously 
$$fr(s, k^+) \quad = \quad fr (s,k)+ fr(s,k+1) + \cdots + fr(s,w).$$ 

Thus everthing hinges on the numbers $fr(s^+,k)$.  They  can be calculated as in Section 3
provided we have\footnote{As previously mentioned, {\it Dualize and Advance} can be employed, and one should also exploit the fact that the facets of ${\cal F}r(s^+)$ help to find the facets of ${\cal F}r((s+1)^+)$.  The matter needs further investigation.} the facets of all simplicial complexes ${\cal F}r(s^+) \ (1 \leq s \leq m)$. 
To fix ideas, consider the concrete database ${\cal D} = \{T_1, \cdots, T_7\}$ given by:

\begin{tabular}{c|c|c|c|c|c|c|c|c|c|}
& 1 & 2 & 3& 4 & 5 & 6 &7 & 8 & 9 \\ \hline 
$T_1$ & $x$ & ${\bf x}$ & $x$ & & $x$ & & ${\bf x}$ & & ${\bf x}$ \\ \hline 
$T_2$ & $x$ & $x$ & $x$ & $x$ & & $x$ & & $x$ & $x$ \\ \hline
$T_3$ & & &  $x$ & & $x$ & $x$ & & $x$ & $x$ \\ \hline
$T_4$ & & ${\bf x}$ & & & $x$ & & ${\bf x}$ & & ${\bf x}$ \\ \hline 
$T_5$ & & & $x$ & & &  $x$ & & $x$ &  \\ \hline 
$T_6$ & & ${\bf x}$ & & $x$ & & & ${\bf x}$ & $x$ & ${\bf x}$ \\ \hline 
$T_7$ & & & $x$ & & & $x$ & & $x$ & $x$ \\ \hline \end{tabular}

Table 6

One calculates these associated numbers $fr(s,k)$ for $0 \leq s \leq 7$ and $1 \leq k \leq 9$:

\begin{tabular}{l|c|c|c|c|c|c|c|c|c|c} 
& 1 & 2& 3 & 4 & 5 &6 & 7 & 8 & 9 & \\ \hline
0 & 0 & 2 & 23 & 69 & 97 & 76 & 35 & 9 &1 & 312\\ \hline
1 & 0 & 13 & 44 & 53 & 29 & 8 &  1& 0 & 0 & 148 \\ \hline 
2 & 2 & 11 & 12 & 3 & 0 & 0 & 0 & 0 & 0 & 28 \\ \hline
3 & 2 & 4 & 4 & 1 & 0 & 0 & 0 & 0 & 0 & 11\\ \hline
4 & 2 & 6 & 1 & 0 & 0 & 0 & 0 & 0 & 0 & 9 \\ \hline
5 & 2 & 0 & 0 & 0 & 0 & 0 & 0 & 0  & 0  & 2 \\ \hline 
6 & 1 & 0 & 0 & 0 & 0 & 0 & 0 & 0 & 0 & 1 \\ \hline
7 & 0 & 0 & 0 & 0 & 0 & 0 & 0 & 0 & 0 & 0 \\ \hline
 & ${9 \choose 1}$ & ${9 \choose 2}$ & ${9 \choose 3}$ & ${9 \choose 4}$ & ${9 \choose 5}$ & ${9 \choose 6}$ & ${9 \choose 7}$ & ${9 \choose 8}$ & ${9 \choose 9}$ & 511 \\ \hline  \end{tabular}
 
Table 7 
 
 For instance, because of $fr(3,3)=4$ there are exactly four $3$-frequent sets of cardinality 3, one of which is $\{2, 7,9\}$ (indicated boldface in Table 6). Many more probabilities can be calculated from Table 7. Say, the probability that a random 3-element itemset is $2^+$-frequent, is
 $$\frac{12+4+1}{{9\choose 3}}\approx 0.202,$$
 whereas the probability that a random $2^+$-frequent set has $3$ elements, is
 $$\frac{12+4+1}{28+11+9+2+1} \approx 0.333.$$
 Similarly one reads from Table 7 that the probability of a $2^+$-frequent set $X$ to be $2$-frequent, is $\frac{28}{51}= 0.549$. If one additionally requires that $|X| \geq 2$ or $X \geq 3$ the corresponding probabilities obviously increase, in fact they are $0.619$ and $0.714$.
 
 {\bf 8.2} \quad Here $s$ will be fixed and we consider the simplicial complex ${\cal S}{\cal C} : = {\cal F}r(s^+)$ with respect to some data base ${\cal D}$. For $X \in {\cal S}{\cal C}$ put $\mbox{supp} (X): = \{T \in {\cal D}: X \subseteq T\}$. One calls $Y \in {\cal S}{\cal C}$ {\it closed} if it is maximal with respect to its support:
 $$(\forall X \in {\cal S}{\cal C}) \quad X \supsetneqq Y \ \Ra \ \mbox{supp} (X) \subsetneqq \ \mbox{supp}\, (Y).$$
 The closed frequent sets play a prominent r\^{o}le in FSM. One verifies at once that they coincide with the members (apart from $W$) of the closure system ${\cal C}$ generated by the facets $F_1, \cdots, F_h$ of ${\cal S}{\cal C}$. 
 
 Given ${\cal D}$, the main result of [UAUA] states that the closed $s^+$-frequent sets can be enumerated with polynomial delay. Here is a quick argument. The closure operator $cl: \ {\cal P}(W) \ra {\cal P}(W)$ coupled to ${\cal C}$ is given by
 $$ cl (X) = \left\{ \begin{array}{lll} \bigcap \ \mbox{supp} (X), & \mbox{if} \ |\mbox{supp} (X)| \geq s\\
 W, & \mbox{otherwise}. \end{array}\right. $$
 Calculating $cl (X)$ costs $O(|{\cal D}|w)$, and so the [GR] method mentioned in Section 4 enumerates the members of ${\cal C}$ with polynomial delay $|{\cal D}|w^2$.
 
From ${\cal C}$ one can get a partition of ${\cal S}{\cal C}$ different from the two partitions considered in 8.1. Namely, ${\cal S}{\cal C}$ is the disjoint union of the set families
  $${\cal S}{\cal C}[Y]: = \{X \in {\cal S}{\cal C}: \ \mbox{supp}(X) = \ \mbox{supp}(Y) \}$$
 when $Y$ runs through ${\cal C}$ (here ${\cal C}$ must be of moderate size).  If $Y_1, \cdots, Y_m$ are all lower covers of $Y$ in ${\cal C}$ then
 $${\cal S}{\cal C}[Y] = \{X \subseteq Y : \ X \not\subseteq Y_1 \ \mbox{and} \ X \not\subseteq Y_2 \ \mbox{and} \, \cdots X \not\subseteq Y_m\}.$$
 Hence starting with the $\{0,2\}$-valued row ${\cal P}(Y)$ and imposing on it the constraints $X \cap (W \backslash Y_i) \neq \emptyset \  (1 \leq i \leq m)$ will deliver ${\cal S}{\cal C}[Y]$ as a disjoint union of few $\{0,1,2,e\}$-valued rows. Such a representation may be useful for certain FSM tasks. Of course it all works for any closure systems ${\cal C}$ (possibly unrelated to FSM). For instance [W3, Section 4.3] is an application geared at speeding up inclusion-exclusion.

\section*{References}

\begin{enumerate}
\item[{[EMG]}] T. Eiter, K. Makino, G. Gottlob, Computational aspects of monotone dualization: A brief survey, Disc. Appl. Math. 156 (2008) 2035-2049.
\item[{[FG]}] I. Flum, M. Grohe, Parametrized complexity theory, Springer-Verlag Berlin Heidelberg 2006.
\item[{[GR]}] B. Ganter, K. Reuter, Finding all closed sets: A general approach, Order 8 (1991) 283-290.
\item[{[K]}] D. Knuth, The Art of Computer Programming, Vol 4, Fascicle 1, Binary Decision Diagrams, Addison-Wesley 2005.
\item[{[KP]}] V. Kaibel, M.E. Pfetsch, Computing the face lattice of a polytype from its vertex-facet incidences, Computational Geometry 23 (2002) 281-290.
\item[{[MS]}] E. Miller, B. Sturmfels, Combinatorial Commutative Algebra, Springer 2005.
\item[{[S]}] R. Stanley, Combinatorics and commutative algebra (2nd ed.), Progress in Math. Vol 41, Birkh\"{a}user 1996.
\item[{[UAUA]}] T. Uno, T. Asai, Y. Uchida, H. Arimura, LCM: An efficient algorithm for enumerating frequent closed item sets, Proc. of the workshop on frequent itemset mining implementations, 2003.
\item[{[W1]}] M. Wild, Counting or producing all fixed cardinality transversals, Algorithmica, DOI 10.1007/s00453-012-9716-5.
\item[{[W2]}] M. Wild, Compactly generating all satisfying truth assignments of a Horn formula, Journal on Satisfiability, Boolean Modeling and Computation 8 (2012) 63-82.
\item[{[W3]}] M. Wild, Inclusion-exclusion meets exclusion, arXiv:1309.6927.
\item[{[W4]}] M. Wild, Output-polynomial enumeration of all fixed-cardinality ideals of poset, respectively all fixed-cardinality subtrees of a tree, Order online April 2013.
\item[{[WK]}] X. Wu, V. Kumar (editors), The Top Ten Algorithms in Data Mininig, Chapman L Hall/CRC 2009.
\end{enumerate}

\end{document}